\documentclass[12pt]{amsart}
\usepackage{pgfplots}
\usepackage{float}
\usepackage[utf8]{inputenc}
\usepackage{graphicx}
\usepackage{color}
\usepackage{mdwlist}
\usepackage{amssymb}
\usepackage{subfigure}
\usepackage{relsize}
\usepackage{graphics}
\allowdisplaybreaks

\usepackage{cleveref}

\newcommand{\beq}{\begin{eqnarray*}}
\newcommand{\feq}{\end{eqnarray*}}
\newcommand{\beqn}{\begin{eqnarray}}
\newcommand{\feqn}{\end{eqnarray}}

\newcommand{\RN}[1]{%
  \textup{\uppercase\expandafter{\romannumeral#1}}%
}

\newtheorem{theorem}{Theorem}[section]

\theoremstyle{definition}

\theoremstyle{remark}
\newtheorem{remark}[theorem]{Remark}
\numberwithin{equation}{section}

\newcommand\numberthis{\addtocounter{equation}{1}\tag{\theequation}}

\begin{document}
\title[ Euler-alignment system]
{Global dynamics of the Euler-alignment system with weakly singular kernel}

\author{Manas Bhatnagar and Hailiang Liu}
\address{Department of Mathematics, Iowa State University, Ames, Iowa 50011}
\email{manasb@iastate.edu}
\email{hliu@iastate.edu} 
\begin{abstract} 
This letter studies the Euler-alignment system with weakly singular influence functions by introducing a novel technique to bound the density. Instead of resorting to a nonlinear maximum principle used in [C. Tan, Nonlinearity, 33: 1907--1924, 2020] to bound the interaction term $\psi\ast\rho$ by $\rho^s$ with $s\in (0, 1)$ for $\psi$ with an algebraic singularity at origin, we bound $\psi*\rho$ by a relaxed constant for any $\psi \in L^1$.  We thus establish the global-in-time existence results with weaker assumptions on $\psi$ and refined solution bounds,  characterized by the structure of $\psi$. 
\end{abstract}
\maketitle

\section{Introduction}
\label{intro}
The Euler-alignment system 
\begin{subequations}
\label{mainsys}
\begin{align}
& \rho_t + (\rho u)_x = 0, \qquad t>0, x\in X,  \label{mainsysmass}\\
& u_t + u u_x = \psi\ast(\rho u) -  u\psi\ast\rho, \label{mainsysmom}
\end{align}
\end{subequations}
as a hydrodynamic model characterizes the self-organized collective bahaviors, in particular alignment and flocking \cite{HL09, HaTa08}, originating from the Cucker-Smale agent-based models \cite{CS07a, CS07b}. Here, $\rho$ represents the density of the group, and $u$ is the associated velocity. 
The term that appears on the right of (\ref{mainsys}b) is the {\it alignment force}, where we use the notation 
$$
\psi\ast \rho = \int_X \psi(|x-y|)\rho(t, y)dy. 
$$
The spatial domain $X$ can be all of $\mathbb R$ or a periodic domain and the initial data is $(\rho(0,x)\geq 0, u(0,x))$. Since $M:=\int_X \rho dx $ is conserved and a simple maximum principle on $u$ indicates that $\int_{X} \rho u$ is well-defined, we have that the right hand side of \eqref{mainsysmom} is well-defined.

For system (\ref{mainsys}) the global dynamics has been  studied for all three scenarios: bounded Lipschitz interaction \cite{CCTT16, TT14}, the strongly singular interaction \cite{DKRT17, KT18}, and the weakly singular interaction \cite{Tan20}. In this letter we make an attempt to extend the results in \cite{Tan20} on the weakly singular interaction, with focus on a new technique circumventing the use of the nonlinear maximum principle in \cite{Tan20}, and also recovering the result in \cite{CCTT16}. 

In order to compare the results in \cite{Tan20}, we recall some conventions here. A sharp critical threshold condition is obtained in \cite{CCTT16} for bounded interaction with the help of an important quantity 
$$
G:=u_x+ \psi\ast \rho, 
$$
the dynamics of which then becomes 
$$
\partial_t G +\partial_x(Gu)=0.
$$
This together with the continuity equation for $\rho$ can serve as an alternative representation of (\ref{mainsys}). The velocity field $u$ can be recovered by the relation for $G$ and the conservation of momentum. We rewrite this new equivalent system here,
\begin{subequations}
\label{eqmainsys}
\begin{align}
& \rho_t + (\rho u)_x = 0, \qquad t>0,  \label{eqmainsysa}\\
& G_t + (u G)_x = 0, \label{eqmainsysb}\\
& u_x = G - \psi\ast\rho. \label{eqmainsysc}
\end{align}
\end{subequations}
For bounded Lipschtz interaction, the critical threshold results in \cite{CCTT16} may be stated as: 
{\it system (\ref{eqmainsys}) admits global regular solution if and only if 
$
\inf_x G(0, x)\geq 0.
$}
In sharp contrast,  the so-called {\it strongly  singular interaction}, say $\psi(r)$ is not integrable at $r=0$, has been shown to have a regularization effect in the sense that global-in-time regular solution is always ensured, see \cite{DKRT17, KT18}. This is true for Euler-alignment, see \cite{DKRT17}, as well as Euler-Poisson-alignment, see \cite{KT18}, systems. For the {\it weakly singular interaction}, the critical thresholds obtained in \cite{Tan20} are about similar to those for the bounded interaction, while the principal global existence states the following.
\begin{theorem}
\label{tan20} 
\cite{Tan20} Consider \eqref{eqmainsys} with smooth initial data $(\rho(0,x)\geq 0, G(0,x))$ and weakly singular interaction $\psi\in L^1(X)$, nonnegative and satisfying, 
$$
\lambda r^{-s} \leq \psi(r) \leq \Lambda r^{-s}, \Lambda \geq \lambda >0,\; s\in (0, 1)
$$
uniformly 
in a neighborhood of the origin. If 
$$
\inf_x G(0, x)>0,
$$
then there exists a globally regular solution to the system \eqref{eqmainsys}. 
\end{theorem}
In this letter, we extend that to the following. 
\begin{theorem}
\label{mainth}
Let $\psi\in L^1(X)$ be non-negative. Consider \eqref{eqmainsys} with smooth initial data $(\rho(0,x)\geq 0, G(0,x))$. If 
$$
\inf_x G(0, x)>0,
$$
then there exists globally regular solution to the system \eqref{eqmainsys}. \\
Moreover, for this solution we have the global bound, 
\begin{align*}
& ||\rho(t,\cdot )||_\infty \leq  \max\left\{ ||\rho(0,\cdot)||_\infty ,\ \beta   \right\},\\
& ||G(t,\cdot )||_\infty  \leq \max\left\{ ||G(0,\cdot)||_\infty,\ \gamma  \right\},
\end{align*}
where $
C_0: = \inf_x\frac{G(0,x)}{\rho(0,x)}>0, 
$
and 
\begin{align*} 
\beta & = \left\{  
\begin{array}{ll}
0 & \text{if} \; C_0>\|\psi\|_{L^1},\\
\frac{M\|\psi\|_\infty}{C_0} & \text{if} \; \psi \in L^\infty,\\
\inf_k \left\{\frac{Mk}{C_0-\int_{\psi \geq k}\psi dx}:\ C_0> \int\displaylimits_{\psi \geq k}\psi dx\right\} & \text{otherwise}, 
\end{array}
\right.
\\
\gamma & = \left\{  
\begin{array}{ll}
M\|\psi\|_\infty & \text{if} \; \psi \in L^\infty(X),\\
\|\psi\|_1  \max\left\{ ||\rho(0,\cdot)||_\infty ,\ \beta   \right\} & \text{otherwise}. 
\end{array}
\right.
\end{align*} 
\end{theorem}
Some remarks are in order.
\begin{remark} For the finite time blow-up, both bounded interaction and weakly singular interaction feature the similar behavior: If $\inf_x G(0, x)<0$, then the solution admits a finite time blow up \cite{CCTT16, Tan20}. For the global-in-time existence, our method unifies the two scenarios as well. 
More precisely, we recover the global regular solution result in \cite{CCTT16}.
\end{remark}
\begin{remark}
A remarkable feature of our approach is that apart from the assumption that $\psi\in L^1$, we only require the nonnegativity of $\psi$. Our analysis neither requires the need of an explicit singularity condition \eqref{l1cond} nor any restriction on the number or position of singularities, and hence, is applicable on the larger class of $L^1$ as well as bounded influence functions. In contrast, the main technique used by the author in \cite{Tan20} is to derive a nonlinear maximum principle, which requires that $\psi$ only has the one singularity at the origin and obeys \eqref{l1cond}.
\end{remark}

The main contribution of this letter is two fold: (i) A novel technique to bound the interaction term is introduced, leading to a global control of the density; (ii) we obtain refined solution bounds and characterize how they are determined by the structure of $\psi$.  Here we not only give an alternative proof to 
the existence of global-in-time solutions, but also do so for more general influence functions. In addition, we are able to establish sharp solution bounds which indicate the effect from the profile of $\psi \in L^1(X)$.   

In Section \ref{Tananalysis}, we review the existing method, and the proof of Theorem \ref{mainth} is presented in Section \ref{analysis}. In Section \ref{example}, we give a particular example of our general result to showcase our technique and the fineness of the bounds obtained. 

\section{Review of the existing method \cite{Tan20}}  
\label{Tananalysis}
It is known from \cite{Tan20} that the local-in-time classical solutions to \eqref{eqmainsys}, and in turn \eqref{mainsys}, can be extended for further times if the following holds,
\begin{align}
\label{bkm}
& \int_0^T ||G(t,\cdot)||_\infty + ||\rho(t,\cdot)||_\infty\, dt < \infty.
\end{align}
Hence, the persistence of smooth solutions upto any time $T>0$ is guaranteed if there are a priori bounds on the infinity norms of $G(T,\cdot),\rho(T,\cdot)$. 

Writing \eqref{eqmainsys} as ODEs along the characteristic path,
\begin{align}
\label{chpath}
& \Gamma = \left\{ (t,\mathcal{X}(t)): \frac{d\mathcal{X}}{dt} = u(t,\mathcal{X}(t)),\ \mathcal{X} (0) = x_0 \right\},
\end{align}
we have,
\begin{subequations}
\label{sysalchpath}
\begin{align}
& \frac{d}{dt}\rho(t,\mathcal X (t)) = -\rho(G-\psi\ast\rho),  \label{sysalchpath1} \\
& \frac{d}{dt} G(t,\mathcal X(t)) = -G(G-\psi\ast\rho), \label{sysalchpath2}
\end{align}
\end{subequations}
with initial point $(\rho_0, G_0):=(\rho(0, x_0), G(0, x_0))$. If $\rho_0=0$, then $\rho(t,\mathcal X(t))\equiv 0$ along the characteristic. Observing \eqref{sysalchpath2}, the bounds on $G$ depend on the global behavior of $\rho$ rather than along a single path. If one is able to bound $\rho,G$ along characteristic paths for which $\rho_0>0$, then the solutions are indeed bounded along all other characteristic paths as well. 
For $\rho_0>0$, we have,
$$
G=\frac{\rho}{\rho_0}G_0,
$$
since  along $\Gamma$ one has, 
$$
\frac{d}{dt}\left( \frac{G(t,\mathcal X(t))}{\rho(t,\mathcal X(t))}\right) = 0. 
$$
Plugging this in \eqref{sysalchpath1}, \begin{equation}\label{rho}
\frac{d}{dt}\rho(t,\mathcal X (t)) = -\frac{G_0}{\rho_0}\rho^2 + \rho\psi\ast\rho.
\end{equation} 
From this we see that if $G_0<0$, then $\rho$ will blow up at finite time since $\frac{d}{dt}\rho(t,\mathcal X (t))>\frac{|G_0|}{\rho_0}\rho^2$, also $G\to -\infty$ at a finite time. For $G_0>0$, from  \eqref{sysalchpath2} we see that 
\begin{align}\label{G}
0\leq G \leq \max \{G_0, \gamma\}
\end{align}
for all $t>0$, if $\psi\ast\rho \leq \gamma$ for some $\gamma>0$. 

It is left to control $\rho$ in the case $G_0>0$ and the problem boils down to ensuring a bound on the term $\psi\ast\rho$. A bounded kernel would immediately result in the needful since the mass is conserved in time. However, a singular kernel could aggravate this term. The author in \cite{Tan20} bounded this term by proving a nonlinear maximum principle on $\rho$. It states that if 
\begin{align}
\label{l1cond}
& \lambda |x|^{-s}\leq\psi(x)\leq \Lambda|x|^{-s},
\end{align}
close to the origin, then for a function $f\in L^1$
$$
(\psi\ast f)(x_\ast) \leq C_1 f^s(x_\ast),
$$
where $s\in(0,1)$ and $x_\ast$ is any point where $f$ attains maximum value.
Evaluating $ \partial_t \rho + u\partial_x \rho=-\rho (G-\psi*\rho)$ at $(t, x^*)$ and using the above result, one obtains 
$$
\partial_t \rho(t,x_\ast) \leq - C_0 \rho^2(t,x_\ast) + C_1 \rho^{1+s}(t,x_\ast),
$$
where $C_0= \inf \frac{G(0,x)}{\rho(0,x)}>0$.  Hence, $\rho$ cannot become very large resulting in being upper-bounded. We can further use this bound in \eqref{G} to obtain that $G$ remains bounded from above if $G_0>0$, thereby obtaining a priori bounds on $||\rho(t,\cdot)||_\infty, ||G(t,\cdot)||_\infty$ for all time. The condition \eqref{bkm} then guarantees global-in-time classical solutions. 

\section{Proof of Theorem \ref{mainth}}
\label{analysis}
\textbf{Step 1:} (Global bound on $\rho$) We use system \eqref{sysalchpath} along the characteristic path \eqref{chpath} and some preliminary notations from Section \ref{Tananalysis} in this Section. Recall $C_0= \inf \frac{G(0,x)}{\rho(0,x)}$ which is strictly positive from the hypothesis of the Theorem. 

For a number $k\geq 0$,
\begin{align}
\label{setB}
& B=\left\{x: \quad \psi(x) \geq k  \right\},
\end{align}
and $B^c$ be its complement on which $\psi \leq k$. Also, a set $B$ is admissible if and only if $4\int_B \psi(x) dx < C_0$. Since $\psi\in L^1$, there exists a $k\geq 0$ for which the set $B$ is admissible. In particular, for all sufficiently large $k$, the corresponding sets $B$ are admissible. However, note that the set $B$ need not be admissible for an arbitrary $k\geq 0$. In particular, if $C_0<||\psi||_1$, then for $k=0$, we have $B=X$ and $C_0\ngtr 4\int_B\psi(x)dx$. 

We can obtain an implicit bound on $\psi\ast\rho$,  
\begin{align*}
\psi\ast\rho & = \int_B \psi(y)\rho(t,x-y)dy +\int_{B^c} \psi(y)\rho(t,x-y)dy\\
& \leq  \|\rho(t, \cdot)\|_\infty \int_B \psi(x) dx  + k \int \rho(t, x)dx\\
& \leq  \|\rho(t, \cdot)\|_\infty \int_B \psi(x) dx  + M k .
\label{psirhoub} \numberthis
\end{align*} 
We shall identify a uniform bound for $\rho$, which is valid for all time $t>0$. Suppose $\|\rho(0, \cdot)\|_\infty \leq a$, then there exists $T=T(a)$ such that $\|\rho(t, \cdot)\|_\infty \leq 2a $ for $t\in [0, T]$. For any $k$ for which the corresponding set $B$ is admissible and for $t\in [0, T]$, we use \eqref{psirhoub} to get,
\begin{align*}
\psi\ast\rho & \leq  2a\int_B\psi(x)dx +Mk  \leq \frac{aC_0}{2} + Mk =: b.
\end{align*}
On the other hand, equation (\ref{rho}) for $\rho$ leads to the following differential inequality 
$$
\frac{d}{dt}\rho(t,\mathcal X (t)) \leq -C_0\rho^2+b\rho, 
$$
for $t\in [0, T]$.  This inequality ensures that along the path $\Gamma$ as in \eqref{chpath},
$$
\rho(t,\mathcal X (t)) \leq \max\left\{ \rho_0, \frac{b}{C_0}\right\}, \quad t\in [0, T]. 
$$
Thus, $\rho(t,\mathcal X (t)) \leq a$  as long as $a \geq \frac{b}{C_0}$.
Therefore, it suffices to take $a$ so that $a \geq   \frac{2Mk}{C_0 }$.
By induction we have, 
\begin{align}
\label{rhoroughbound}
& \rho(t,\mathcal X(t)) \leq \max \left\{\rho_0, \frac{2Mk}{C_0}\right\}.
\end{align}
\textbf{Step 2:} (Refined bound on $\rho$) After we know $\rho$ is globally bounded, we can further improve the bound in two steps. Again assume this bound to be $a$, then we have 
$$
\psi\ast\rho \leq a \int_A \psi(x) dx  + M k.
$$
For $k\geq 0$, let
\begin{align}
\label{setA}
A = \left\{x: \quad \psi(x) \geq k  \right\}.
\end{align}
Also, a set $A$ is admissible if and only if $\int_A \psi(x) dx < C_0$.
Note that by relaxing the admissibility condition on $A$ compared to $B$ in \eqref{setB}, we allow for more, and in particular, smaller values of $k$, which will result in an even sharper bound.
We have, 
$$
\frac{d}{dt}\rho(t,\mathcal X (t)) \leq -C_0\rho^2+\rho \left(a \int_A \psi(x) dx  + M k\right), 
$$
where $a$ is a bound on $\rho(t,\mathcal X(t))$ which we know exists from Step 1. Consequently,
\begin{align*}
  \frac{d}{dt}\rho(t,\mathcal X (t)) & \leq -C_0\rho \left( \rho - \frac{a \int_A \psi(x) dx  + M k}{C_0}\right) 
\end{align*}
Therefore, for any $k$ for which the set $A$ is admissible, it suffices to take an $a$ such that,
$$
a \geq  \max\left\{\rho_0, \frac{a \int_A \psi(x) dx  + M k}{C_0}\right\}.
$$
Hence,
$$
a\geq  \max\left\{ \rho_0, \frac{ M k}{C_0-\int_A\psi(x)dx}\right\}.
$$
Consequently, we have a global bound on density,
\begin{align}
\label{rhoroughbound2}
& \|\rho(t,\cdot)\|_\infty \leq \max\left\{ \|\rho(0,\cdot)\|_\infty,  \frac{ M k}{C_0-\int_A\psi(x)dx}\right\} .
\end{align}
For fixed $k$, \eqref{rhoroughbound2} is a better bound than \eqref{rhoroughbound}. 

\textbf{Step 3:} (Further narrowing down by minimization in $k$) Next we examine how the profile of $\psi$ can affect this bound, and how to even further sharpen the bound. Firstly, observe that if $C_0>\|\psi\|_1$, then $k$ can be zero with $A$ as the complete space being an admissible set from Step 2. This results in a maximum principle on $\rho$,
$$
\|\rho(t,\cdot)\|_\infty \leq \|\rho(0,\cdot)\|_\infty.
$$
Also, if $\psi$ is bounded, then for any $\epsilon >0$ we can have $k=\max\psi+\epsilon$ with $A=\emptyset$. By taking limit as $\epsilon\to 0$, we obtain,
$$
\|\rho(t,\cdot)\|_\infty \leq \max\left\{ \|\rho(0,\cdot)\|_\infty,  \frac{ M \max\psi}{C_0}\right\}.
$$
This points to the fact that the quantity $\frac{ M k}{C_0-\int_A\psi(x)dx}$ adjusts itself as and when $\psi$ changes. This motivates us to treat the following optimization problem to further sharpen the bounds on density and see how $\psi$ affects these bounds,
\begin{align}
\label{infprob}
& \beta  : = \inf\left\{ \frac{ M k}{C_0-\int_{A_k}\psi(x)dx}:\ C_0> \int_{A_k}\psi(x)dx  \right\}.
\end{align}
Note that we now explicitly show the dependence of the set $A$ on $k$ by the subscript. We consider the case where $C_0\leq \|\psi\|_1$ for otherwise, $\beta  = 0 $ (minimum in this case) trivially.
We set,
$$
g(k):= \frac{ M k}{C_0-\int_{A_k}\psi(x)dx}, 
$$
with domain such that $C_0 > \int_{A_k}\psi(x)dx$. In particular, let 
$$
k_0 := \sup \left\{ k: C_0 = \int_{A_k}\psi(x)dx \right\} \geq 0 .
$$ 
Then $g:(k_0,\infty )\to\mathbb{R}^+$ is a well-defined function. Also, $\lim_{k\to\infty} g(k) = \infty$. This combined with the fact that $g$ is bounded from below, we have that \eqref{infprob} is finite and the minimizer if exists is contained in the set $[k_0,\infty)$. Consider the sequence $\{k_n\}_{n=1}^\infty$ such that $g(k_n)\searrow \beta $.
We have,
$$
||\rho (t,\cdot)||_\infty \leq \max \left\{||\rho(0,\cdot)||_\infty,\ g(k_n) \right\},\quad n\geq 1.
$$
Taking limit as $n\to\infty$, we obtain our final optimized bound,
$$
||\rho (t,\cdot)||_\infty \leq \max \left\{||\rho(0,\cdot)||_\infty, \beta \right\},\quad t>0.
$$
With this bound we have 
$$
\psi\ast \rho \leq ||\rho (t,\cdot)||_\infty \|\psi\|_1 \leq 
\max \left\{||\rho(0,\cdot)||_\infty, \beta \right\} \|\psi\|_1=:\gamma. 
$$
This when combined with (\ref{G}) gives the bound for $G$. Note that for $\psi\in L^\infty(X)$, we simply take $\gamma=M \|\psi\|_\infty.$ In summary,  
the uniform bounds for all time $t>0$ read as,
\begin{subequations}
\label{rhoGub}
\begin{align}
& ||\rho (t,\cdot)||_\infty \leq \max \left\{||\rho(0,\cdot)||_\infty, \beta  \right\}, \label{rhoGub1}\\
& ||G(t,\cdot)||_\infty \leq \max\left\{ ||G(0,\cdot)||_\infty,\ \gamma  \right\}, \label{rhoGub2}
\end{align}
\end{subequations}
where 
$$
\gamma= \left\{  
\begin{array}{ll}
M\|\psi\|_\infty & \text{if} \; \psi \in L^\infty(X),\\
\|\psi\|_1 \max \left\{||\rho(0,\cdot)||_\infty, \beta  \right\} & \text{otherwise}.
\end{array}
\right.
$$
This finishes the proof of Theorem \ref{mainth}. \qed \\

\section{An example}
\label{example}
In this section, we show the importance of optimization with level set $\psi=k$ by explicitly calculating $\beta $ for some specific choices of alignment kernels with finite weight. It allows us to get an insight into the optimized bounds on density.

We consider an example: $\psi=|x|^{-\alpha}$ with $\alpha\in (0,1)$ on the periodic domain $X = \mathbb{T} = (-1/2,1/2]$,  hence $\|\psi\|_1=\frac{2^\alpha}{1-\alpha}$. We have two cases relative to the weight of $\psi$.

\textbf{Case 1:} If $C_0> \frac{2^\alpha}{1-\alpha}$, then corresponding to $k=0$, we have $A=\mathbb{T}$ is admissible. Therefore, the minimizer is $k=0$ and consequently, $\beta=0$. We then have the maximum principle on $\rho$.

\textbf{Case 2:} If $0 < C_0\leq \frac{2^\alpha}{1-\alpha}$, then we need to minimize the function,
$$
\beta(k) = \frac{Mk}{C_0 - \int_{A_k}\psi(x) dx},
$$
with the domain being all $k\in\mathbb{R}^+$ such that $C_0>\int_{A_k}\psi(x)dx$. Note that 
$$
A_k = \{ x: |x|^{-\alpha} \geq k\}=\{x, |x|\leq k^{-1/\alpha} \}.
$$
Therefore,
\begin{align*}
\int_{A_k}\psi(x) dx & = 2\int_0^{k^{-1/\alpha}} x^{-\alpha} dx = \frac{2k^{1-1/\alpha}}{1-\alpha}.
\end{align*}
Let $k_0$ be such that $\int_{A_{k_0}}\psi(x) dx=C_0$, then the minimization problem reduces to the following: 
\begin{align*}
& \min_{\left( k_0 ,\infty \right)} \beta(k) := \min_{\left( k_0 ,\infty \right)} \frac{Mk}{C_0 - \frac{2k^{1-1/\alpha}}{1-\alpha}},
\end{align*}
where $k_0 = \left(\frac{2}{C_0(1-\alpha)}\right)^{\frac{\alpha}{1-\alpha}}$.
Let the minimizer of $\beta(k)$ be $k^\ast$. We have $\frac{d\beta}{d k} (k^\ast) = 0$. Hence,
$$
k^\ast = \left(\frac{2}{C_0 \alpha(1-\alpha)}\right)^{\frac{\alpha}{1-\alpha}}.
$$
Note that since $\alpha\in (0,1)$, we have $k^\ast > k_0$. Consequently, the minimizer is unique. Consequently,
$$
\beta:= \beta(k^\ast) = \left( \frac{2}{\alpha} \right)^{\frac{\alpha}{1-\alpha}}\frac{M}{(C_0(1-\alpha))^{\frac{1}{1-\alpha}}}.
$$
In the special case when $\alpha = 1/2$, we have $\beta=16M/C_0^2$. Hence the bound on density would be,
$$
||\rho(t,\cdot)||_\infty \leq \max\left\{ ||\rho(0,\cdot)||_\infty, \frac{16M}{C_0^2} \right\}, \quad t>0.
$$

\section*{Acknowledgments}
This work was supported by the National Science Foundation under Grant DMS1812666.

\bigskip

\bibliographystyle{abbrv}

\end{document}